\documentclass{article}
\usepackage{amssymb, amsmath, amsthm, bbm, verbatim}
\usepackage{color}

\newtheorem{teo}{Theorem}

\renewcommand{\r}{{\mathbb R}}

\title{On alternative wavelet reconstruction formula: a case study of approximate wavelets
}
\author{
Elena A.  Lebedeva
\footnote{Mathematics and Mechanics Faculty, St. Petersburg State University,
Universitetsky prospekt, 28, Peterhof,  Saint Petersburg,
 198504, Russia; Institute of Applied Mathematics and Mechanics, St. Petersburg State Polytechnical University,
 Polytechnicheskay 29, 195251, St. Petersburg, Russia},
 Eugene B. Postnikov
\footnote{Department of Theoretical Physics,
 Kursk State University, Radishcheva st. 33, Kursk 305000, Russia}
}
\date{
 ealebedeva2004@gmail.com, postnicov@gmail.com}

\begin{document}

\maketitle

\begin{abstract}
The application of the continuous wavelet transform to study of a wide class of physical processes with oscillatory dynamics is restricted by large central frequencies due to the admissibility condition. We propose an alternative reconstruction formula for the continuous wavelet transform, which is applicable even if the admissibility condition is violated. The case of the transform with the standard Morlet wavelet, which is an important example of such analyzing functions, is discussed.
\end{abstract}



\section{Introduction}

The continuous wavelet transform (CWT) on $L_2(\r)$ 
is defined as 
\begin{equation}
W_{n,\psi} f(a,\,b)=\int_{\mathbb{R}} f(x) \overline{\psi_{n,a,b}(x)}\,dx,
\label{wvlT}
\end{equation}
where 
$n=1$ or $n=2,$ 
$$
\psi_{1,a,b}(x)=\frac{1}{|a|}\psi\left(\frac{x-b}{a}\right), \quad
\psi_{2,a,b}(x)=\frac{1}{|a|^{1/2}}\psi\left(\frac{x-b}{a}\right),
$$ 
$a, b \in\r,$ $a\neq 0.$
For $n=1$, we suppose $\psi \in L_1(\r)\cap L_2(\r)$ and $\int_{\r}|\psi(x)|\,dx=1$,
then the amplitude norm  $\int_{\r}|\psi_{1,a,b}(x)|\,dx=1.$
For $n=2$, we suppose $\psi \in  L_2(\r)$ and $\int_{\r}|\psi(x)|^2\,dx=1$,
then the energy norm $\left(\int_{\r}|\psi_{2,a,b}(x)|^2\,dx\right)^{1/2}=1.$

The CWT is one of powerful modern analysis tools in various branches of science connected with the processing of non-stationary signals since it allows to obtain a detailed localized time-frequency decomposition of non-stationary signals, see e.g.\cite{Addison2002book} and references therein.

At the same time, CWT can be applied not only for an analysis and for the decomposition of signals but also for their reconstruction using time-localized oscillating components. In particular, this framework finds actual implementations in the modern problems of quantum mechanics since wavelets provide a natural way to represent coherent  states corresponding to the mentioned wavelet-based reconstruction \cite{Muzhikyan2012,Toutounji2013}, and neurodynamics \cite{Pavlov2012}, where wavelet-like spikes are typical elements of detected activity.

The opportunity for the wavelet reconstruction is provided by the admissibility condition \cite{Daubechies1992book}
$$
C_{\psi}=\int_{\r}\frac{|\widehat\psi(\omega)|^2}{|\omega|}d\omega<\infty
$$ 
that is equivalent to 
 $\widehat\psi(0)=0$ under the additional assumption 
 $(1+|\cdot|^{\alpha})\psi \in L_1(\r),$ $\alpha >0.$ 
As usual, we denote by $\widehat{f}$ the Fourier transform of $f$: 
$\widehat{f}(\omega)=\int_{\r}f(x)e^{i \omega x}\,dx.$ 

Under the admissibility condition the reconstruction formula takes place for all $f\in L_2(\r)$
 \begin{equation}
\label{classic1}
	f(x)=\frac{1}{C_{\psi}}\int_{\r}\int_{\r} \psi_{2,a,b}(x) W_{2,\psi}f(a,\,b) \frac{dadb}{a^2},
\end{equation}
where the equality is understood in a weak sense. If, in addition, $f$ is continuous at $x\in\r$, then (see \cite[Th.3.10]{Chui}) the equality holds at the point $x\in\r$.  The analogous result takes place for a pair of different wavelets $\psi, g$ used for the analysis (the function $\psi$) and the reconstruction (the function $g$). In this case the admissibility condition is 
$
\omega^{-1}\overline{\widehat{\psi}(\omega)}\widehat{g}(\omega) \in L_1(\r),
$ 
and
$
C_{\psi,g}:= \int_{\r} |\omega|^{-1}\overline{\widehat{\psi}(\omega)}\widehat{g}(\omega)\,d\omega \neq 0.
$
And the result itself is read as follows (see \cite[Prop.2.4.2]{Daubechies1992book}):
if $\psi,$ $g,$ $x g(x) \in L_1(\r),$ $g'\in L_2(\r)$, $\widehat{\psi}(0)=\widehat{g}(0)=0,$  $f\in L_2(\r)$ is bounded, $f$ is continuous at $x\in\r,$   
then  we get
\begin{equation}
\label{classic2}
f(x)=\frac{1}{C_{\psi,g}} \lim_{A_1\to 0,\,A_2\to\infty} \int_{A_1\leq|a|\leq A_2}\frac{da}{a^2} \int_{\r}
g_{2,a,b}(x) W_{2,\psi}f(a,\,b)\,db
\end{equation}
The function can also be reconstructed (see \cite[Th.(2.2)]{Holschneider1991}) as 
\begin{equation}
\label{classic3}
f(x)=\lim_{\rho\to\infty \, \varepsilon\to 0}
\int_{\varepsilon}^{\rho}\frac{da}{a}\int_{\r}g_{1,a,b}(x)W_{1,\psi}f(a,\,b)\,db
\end{equation}
under the following restrictions on the function $f$, and the  wavelets $\psi$,  $g$:

	a) $f$ is continuous at $x\in\r,$ 

	b) $\lim_{u\to\infty}\frac{1}{2u}\int_{t-u}^{t+u}f(x)\,dx=0$ uniformly in $t$;

c)
$\log(2+|\cdot|)\psi \in L_1(\r),$ and $\log(2+|\cdot|)g \in L_1(\r),$

 d)
 $\omega^{-1}\widehat{\psi}(\omega)\overline{\widehat{g}(\omega)}\in L_1(\r),$ 
 
 e)
 $\int_0^{\infty}\omega^{-1}\widehat{\psi}(\omega)\overline{\widehat{g}(\omega)}\,d\omega =\int_{-\infty}^0 |\omega|^{-1}\widehat{\psi}(\omega)\overline{\widehat{g}(\omega)}\,d\omega=1.$

However  items d), e)  are  more restrictive than the conditions $C_{\psi,g}<\infty,$ $C_{\psi,g}\neq 0,$ and especially $C_{\psi}<\infty$ ($\psi=g$) coinciding in this case with d).

Note that the a.e. and $L_p(\r)-$norm ($1 \leq p < \infty$) convergence of the inversion formulae is studied intensively  in \cite{Li2012, Rao1994, Weisz2013}. 

At the same time, even one of the most popular and useful in applications wavelets, the standard Morlet wavelet $\psi_M(\xi)=\exp(i\omega_0)\exp(-\xi^2/2)$ does not satisfy the admissibility condition since $\widehat{\psi_M}(0)= C_1 \exp(-\omega_0^2/2)$. However this quantity is sufficiently small for the practically used central frequencies (usually $\omega_0\geq 5$) that allows to apply it widely to the signal decompositions, when an exact reconstruction is not necessary \cite{DeMoortel2004, Addison2006, Postnikov2007, Postnikov2009}.

We do not discuss the deep celebrated background of the classical wavelet reconstruction formula such as the Calderon's reproducing formula, Hilbert spaces with reproducing kernels, irreducible unitary representations of Lie groups.
Our goal is to prove that the violation of the admissibility conditions actually do not forbid an existence of the exact wavelet inversion. As a result, we suggest an alternative wavelet reconstruction formula that does not require
 the admissibility condition. As an example, we consider the standard Morlet wavelet for arbitrary central frequencies including the limiting case $\omega_0\to 0$ when the continuous wavelet transform reduces to the result of diffusion smoothing of a processed signal. We use the Gaussian function to show the difference between the suggested formula and the classical one.

\section{An alternative reconstruction formula}

\begin{teo}
\label{new_inv}
If $f,$ $\psi,$  $\omega \widehat{\psi}(\omega) \in L_2(\r)$ and 
$\widehat{f},$ $\widehat{\psi} \in L_1(\r),$ 
then a. e. on $\r$  
\begin{eqnarray}
	\frac{1}{\pi}{\rm v.p.}\int_{\r}\frac{db}{b-x} \int_{\r}\frac{\partial}{\partial b} W_{1,\psi} f(a,\,b)\,da = \overline{\psi(0)}\, f(x), \nonumber \\
	&& \label{rec1} \\
	\frac{1}{\pi}{\rm v.p.}\int_{\r}\frac{db}{b-x} \int_{\r} \sqrt{|a|} \frac{\partial}{\partial b} W_{2,\psi} f(a,\,b)\,da = \overline{\psi(0)}\, f(x).\nonumber
\end{eqnarray}
 In addition, if 
${\rm supp} \widehat{f}\subset [0,\,\infty),$ then a. e. on $\r$
\begin{equation}
\label{rec2}
	- i \int_{\r}\frac{\partial}{\partial b} W_{1,\psi} f(a,\,b)\,da = \overline{\psi(0)}\, f(b),
	\quad
	- i \int_{\r}\sqrt{|a|}\frac{\partial}{\partial b} W_{2,\psi} f(a,\,b)\,da = \overline{\psi(0)}\, f(b).
\end{equation}
\end{teo}

\noindent {\bf Proof.} We consider the case of the amplitude norm. The case of the energy norm can be proven analogously.  It follows from the definition of CWT that
$$
W_{1,\psi}f (a,\,b)=\int_{\mathbb{R}} f(x) \overline{\psi_{a,b}(x)}\,dx= \frac{1}{2\pi}
\int_{\mathbb{R}} \widehat{f}(\omega) \overline{\widehat{\psi}(a\omega)} e^{i \omega b}\,d\omega
$$
Since $\int_{\mathbb{R}}\left|\widehat{f}(\omega)\overline{\widehat{\psi}(a \omega)} \omega\right|\, d\omega < \infty$ for any fixed $a\in\mathbb{R},$ we get
$$
- i \frac{\partial}{\partial b} W_{1,\psi}f (a,\,b)= \frac{1}{2\pi} \int_{\mathbb{R}}\widehat{f}(\omega)\overline{\widehat{\psi}(a \omega)} \omega e^{i \omega b}\,d\omega.
$$
Using $\widehat{f},$ $\widehat{\psi} \in L_1(\mathbb{R}),$ we obtain
$$
\int_{\mathbb{R}}d\omega \,\int_{\mathbb{R}} \left|\widehat{f}(\omega)\overline{\widehat{\psi}(a \omega)} \omega\right|\, da = \int_{\mathbb{R}} \left|\widehat{f}(\omega)\right|\, d\omega  \int_{\mathbb{R}}\left|\widehat{\psi}(\xi)\right|\, d\xi < \infty.
$$
Whence by the Fubini theorem we finally get
\begin{multline*}
- \int_{\mathbb{R}}\frac{\partial}{\partial b} W_{1,\psi}f(a,\,b)\,da = \frac{- i}{2\pi}
\int_{\mathbb{R}} \widehat{f}(\omega) e^{i \omega b}\,d\omega 
 \int_{\mathbb{R}} \overline{\widehat{\psi}(a \omega)} \omega \,d a \\
 = \frac{\overline{\psi(0)}}{2\pi} \int_{\mathbb{R}} (-i){\rm sgn}(\omega) \widehat{f}(\omega) e^{i \omega b}\,d\omega=
 \frac{\overline{\psi(0)}}{2\pi} \int_{\mathbb{R}} \widehat{Hf}(\omega) e^{i \omega b}\,d\omega=
 \overline{\psi(0)}\, Hf(b), 
\end{multline*}
where $H$ is the Hilbert transform.
It is well known that the Hilbert transform is invertible operator on $L_2(\r),$ $H^{-1}=-H,$ and the inversion formula $H^{-1}(Hf)(x)=f(x)$ holds true almost everywhere on $\r$. Therefore (\ref{rec1}) is proved.

If in addition, ${\rm supp}\widehat{f} \subset[0,\,\infty),$ then the last chain of  equalities is rewritten as follows
$$
- i \int_{\mathbb{R}}\frac{\partial}{\partial b} W_{1,\psi}f(a,\,b)\,da = 
\frac{\overline{\psi(0)}}{2\pi} 
\int_{\mathbb{R}} \widehat{f}(\omega) e^{i \omega b}\,d\omega=
\overline{\psi(0)}\, f(b). 
$$
Therefore (\ref{rec2}) is proven. \hfill $\Box$

To illustrate the result of Theorem \ref{new_inv}, we consider reconstruction formula (\ref{rec2}), choose $n=1,$ and suppose  
$\psi, \, f\in S,$ where $S$ is the Schwartz space. Then $\widehat{f}(\cdot) \overline{\widehat{\psi}(a\cdot)}\in S$, therefore 
$W_{1,\psi}f (a,\,\cdot)\in S$
 for any fixed $a\in \r$, $a\neq 0.$ 
Thus, the derivative $\frac{\partial}{\partial b} W_{1,\psi}f(a,\,b)$ 
is represented using the derivative of the Dirac delta function $\delta$ 
$$
\frac{\partial}{\partial b} W_{1,\psi}f(a,\,b)= 
- \int_{\r}  W_{1,\psi}f(a,\,t) \delta'(t-b)\, dt.
$$
One can identify $\delta'(t-b)$ and $\delta'_{1,a,b}(t)$. 
Indeed, for any $\phi \in S$ we get
$$
\left(\delta_{1,a,b},\,\phi\right)=\int_{\r}\delta_{1,a,b}(x)\phi(x)\,dx=
\int_{\r}\delta(x)\phi(ax+b)\,dx =\phi(b)=\left(\delta(\cdot-b),\,\phi\right),
$$
therefore,
$$
\left(\delta_{1,a,b}',\,\phi\right)= - \left(\delta_{1,a,b},\,\phi'\right) = 
- \left(\delta(\cdot-b),\,\phi'\right) = \left(\delta'(\cdot-b),\,\phi\right).
$$
Thus if ${\rm supp}\,\widehat{f}\subset[0,\,\infty),$ then we obtain "quasi classical" form for the reconstruction formula (\ref{rec2})
\begin{equation}
\label{quasicl}
	f(x)=\frac{i}{\overline{\psi(0)}}\int_{\r}\int_{\r} 
	W_{1,\psi}f(a,\,b) \delta'(b-x)\,db \,da
\end{equation}

However, 
(\ref{quasicl}) is not reduced to the known reconstruction results (\ref{classic1})-(\ref{classic3}). In fact,  admissibility conditions are violated under the choice $\psi(x)=C e^{-c x^2}$ and $g=\delta'$. Indeed $\widehat{\delta'}(\omega)=i\omega$ and 
$\widehat{\psi}(\omega)=C_1 e^{-c_1 \omega^2}$, $c_1>0$
$$
C_{\psi,g} = \int_{\r} |\omega|^{-1}\overline{\widehat{\psi}(\omega)}\widehat{g}(\omega)\,d\omega=
 i C_1 \int_{\r} {\rm sgn}(\omega) e^{-c_1 \omega^2}\, d\omega =0.
$$
The condition e) cited in the Introduction is not fulfilled as well:  
$$
\int_0^{\infty}\omega^{-1}\widehat{\psi}(\omega)\overline{\widehat{g}(\omega)}\,d\omega =-\int_{-\infty}^0 |\omega|^{-1}\widehat{\psi}(\omega)\overline{\widehat{g}(\omega)}\,d\omega
\left(=
i C_1 \int_0^{\infty} e^{-c_1 \omega^2} \, d \omega\right).
$$

On the other hand, 
consider the standard Morlet wavelet with the central frequency 
$\omega_0$, that is $\psi_M(x)= C \exp(i\omega_0 x)\exp(-x^2/2)$, where $C$ is a
 factor depending on chosen norm, and $\widehat{\psi_M}(\omega)= C_1 \exp(-(\omega-\omega_0)^2/2).$  It is clear that $C_{\psi_M}=\infty,$  and 
 the item e) is not fulfilled (set here $g = \psi = \psi_M$). So, the reconstruction formulae (\ref{classic1}) and (\ref{classic3}) are not applicable. 
However, $\psi_M$ satisfies all conditions of Theorem \ref{new_inv}
so one can restore a function $f$ by (\ref{rec1}) and (\ref{rec2}). 
This fact can be interpreted by means of the  "quasi classical" formula (\ref{quasicl}). It is sufficient to set formally $\psi_M$ as the analyzing wavelet and $\delta'$ as  the reconstruction wavelet. In this case $C_{\psi_M, \delta'}<\infty,$
$C_{\psi_M, \delta'}\neq 0$.


\section{Discussion of numerical realization and possible applications}

The proposed method of the inversion has also advantages in a numerical realization. In practice, the result of the transform is represented as a matrix $Wf_{a_j,b_k}$ whose derivative can be easily found by finite difference formulae with a prescribed accuracy. Moreover, the widely used mathematical software, e.g. MATLAB, allows to evaluate this operation in parallel for all columns of the entire matrix. In the same way, the parallel numerical integration  of the obtained result applied for all columns lead to the required reconstruction. 

Since one of most popular fast numerical methods for the continuous wavelet transform utilizes Fast Fourier Transform as an intermediate step \cite{Torrence1998}. Thus both initial sample $f(t_k)$ and transform's matrix $Wf_{a_j,b_k}$ are periodized with respect to the nodes indexed by $k$. Therefore the numerical differentiation  in the form (for example, of an equispaced sample with the step $h=b_{k+1}-b_{k}$, second order of accuracy)
\begin{eqnarray*}
\frac{\partial}{\partial b} Wf_{a_j,b_k}\approx \left(Wf_{a_j,b_{k+1}}-Wf_{a_j,b_{k-1}}\right)/2h,&2<j<k_{max},\\
\frac{\partial}{\partial b} Wf_{a_j,b_1}\approx \left(Wf_{a_j,b_{2}}-Wf_{a_j,b_{k_{max}}}\right)/2h,&\\
\frac{\partial}{\partial b} Wf_{a_j,b_{k_{max}}}\approx \left(Wf_{a_j,b_{1}}-Wf_{a_j,b_{k_{max}-1}}\right)/2h,&\\
\end{eqnarray*}
keeps the periodization property and provides the maximal accuracy within the error of approximation.

The standard Morlet wavelet with a small central frequency, which strictly violates admissibility condition, contains only small number of oscillations inside of the Gaussian envelope. Therefore, it is best adjusted to the extraction of temporal dynamics of emergence and moving localization for short-time pulses, i.g. of acoustic echo \cite{Addison2006} or spike trains in neuroscience \cite{Pavlov2012}

The variation of the central frequency of the Morlet wavelet, which tends to decreasing $\omega_0$ allows not only more fully characterize damping properties of oscillations \cite{DeMoortel2004} but also to reveal topological properties of attractors in the theory of phase synchronization of chaotic oscillators \cite{Postnikov2007, Postnikov2009}.

Moreover, the sufficiently non-restrictive conditions for the function $\psi$ used in the Theorem~\ref{new_inv} indicates that this method of reconstruction is not limited by the wavelets (in a general sense) but may be applied to the appropriate non-oscillating analyzing kernel functions as well, for example, to those Gaussian-based, which emerge in the diffusion signal and image time(space)-frequency smoothing and processing \cite{Gosme2005}.

\section*{Acknowledgments}
The first author is supported by RFBR 12-01-00216-a.
The second author is supported by the grant N0. 1391  of the Ministry of Education and Science of the Russian Federation within the basic part of research funding No. 2014/349 assigned to Kursk State University.

\end{document}